\documentclass{article}

\usepackage{graphicx} 
\usepackage{amsmath,amsfonts,amssymb}
\usepackage{color}

\newtheorem{thm}{Theorem}[section]

\newcommand{\vect}[1]{\boldsymbol{#1}}
\title{Convergence of a meshless numerical method
for a chemotaxis system with density-suppressed motility}
\author{F. Herrero-Hervás, M. Negreanu, A. M. Vargas}
\date{March 2023}

\begin{document}

\maketitle

\begin{abstract}
This article studies a parabolic-elliptic system modelling the pattern formation in E. coli bacteria in response to a chemoattractant known as acylhomoserine lactone concentration (AHL). The system takes into account certain bacterial strains with motility regulation, and the parameters of the equations represent the bacterial logistic growth, AHL diffusion and the rates of production and degradation of AHL. We consider the numerical solution to the system using the Generalized Finite Difference (GFD) Method, a meshless method known to effectively compute numerical solutions to nonlinear problems. The paper is organized to first explain the derivation of the explicit formulae of the method, followed by the study of the convergence of the explicit scheme. Then, several examples over regular and irregular meshes are given.
\end{abstract}

\section{Introduction}\label{intro}
Chemotaxis is a biological process through which certain biological species direct their movements in response to a chemical gradient, either towards the highest or lowest  concentrations of the substance. In these cases, the chemotaxis is said to be positive (and the substance hence called a chemoattractant) or negative (chemorepellent), respectively. Since the seminal works of Keller and Segel \cite{KS1, KS2}, the nonlinearities of the form $-\nabla \cdot (\chi u \nabla v)$ have attracted great interest, as can be found in the surveys by Horstmann, Bellomo et al., Hillen and Painter \cite{Bellomo, Hillen, Horstmann}.
\\\\
The system studied in this article was initially introduced by Liu et al \cite{Liu} to model pattern formation in \textit{E. coli} bacteria in response to a chemoattractant known as acylhomoserine lactone concentration (AHL), which is produced by bacteria themselves. Specifically, certain bacterial strains with motility regulation are considered, so that the system is modeled by the following parabolic equations:
\begin{equation} \label{1.1}
\begin{cases}
    \displaystyle \frac{\partial u}{\partial t} = \Delta (\gamma(v) u) + \mu u \left(1 - \frac{u}{\rho_s}\right), \\\\
    \displaystyle \frac{\partial v}{\partial t} = D_v \Delta v - \alpha v + \beta u,
\end{cases}
\end{equation}
where $u$ represents the bacterial density, $v$ the concentration of AHL and $\mu$, $\rho_s$, $D_v$, $\alpha$ and $\beta$ are positive parameters for the bacterial logistic growth, AHL diffusion and the rates of production and degradation of AHL, respectively. Moreover, $\gamma$ is a known bounded function of $v$ modelling the motility regulation of the considered bacterial strains. 
\\\\
As the diffusion coefficient for AHL is large enough with respect to the rest of the parameters, in Tello \cite{tello}, the second equation of \eqref{1.1} is approximated by an elliptic equation. After rescaling the variables, the following system is obtained:
\begin{equation} \label{1.2}
\begin{cases}
\displaystyle  \frac{\partial u}{\partial t} - \Delta(\gamma(v) u) = \mu u (1-u), \\\\
-\Delta v + v =  u.
\end{cases}
\end{equation}
System \eqref{1.2} is considered over a regular bounded domain $\Omega$, with homogeneous Neumann boundary conditions and appropriate initial data:
\begin{equation} \label{1.3}
\begin{cases}
    u(\vect{x},0) = u_0(\vect{x}), \quad \vect{x} \in \Omega \\\\
    \displaystyle \frac{\partial u}{\partial \vec{n}}(\vect{x},t) = \frac{\partial v}{\partial \vec{n}}(\vect{x},t) = 0, \quad \vect{x} \in \partial \Omega, \hspace{0.3 cm} t > 0.
\end{cases} 
\end{equation}
Futhermore, in \cite{tello}, the author proves that a unique solution globally exists in time, verifying 
\begin{equation} \label{1.4.1}
\lim_{t\to\infty} ||u-1||_{L^\infty(\Omega)} + ||v-1||_{L^\infty(\Omega)} = 0    
\end{equation}
if the following hypothesis are satisfied:
\begin{equation} \label{1.4}
    \begin{cases}
        \gamma \in C^3([0, \infty)), \\\\
       \gamma(s) \geq 0, \hspace{0.2 cm} \gamma'(s) \leq 0, \hspace{0.2 cm} \gamma''(s) \geq 0, \hspace{0.2 cm} \gamma'''(s) \leq 0, \\\\
       -2 \gamma'(s) + \gamma''(s) s \leq \mu_0 < \mu,\\\\
       \displaystyle \frac{|\gamma'(s)|^2}{\gamma(s)} \leq c_\gamma \leq \infty,
    \end{cases} \hspace{0.5 cm} \text{for all } s\geq 0,
\end{equation}
and
\begin{equation}\label{1.5}
\begin{cases}
    u_0 \in C^{2, \alpha}(\overline{\Omega}), \\\\
    \displaystyle \frac{\partial u_0}{\partial \Vec{n}} = 0 \text{ in } \partial \Omega \\\\
    0 < \underline{u}_0 < u_0 < \overline{u}_0 < \infty,
\end{cases}    
\end{equation}
for certain positive constants $\mu_0$, $c_\gamma$, $\overline{u}_0$, $\underline{u}_0$.
\\\\
In our work, we study the numerical solution to system \eqref{1.2}, under conditions \eqref{1.4}-\eqref{1.5}. To do so, we consider the Generalized Finite Difference (GFD) Method, a meshless method widely used to effectively compute numerical solutions to nonlinear problems, many of which are chemotactic processes. For instance, in \cite{bbg}, the authors proved the convergence of the GFD formulae to the periodic solution to a Keller-Segel model with logistic source and negligible diffusion of the species with respect to the substance, and its generalisation in \cite{nv}. Also, the parabolic-parabolic case with periodic enviromental conditions was studied in \cite{bbg1}. The local stability of the constants equilibrium solutions and the GFD numerical solution to a system modelling chemotaxis-haptotaxis was obtained in \cite{bbg2}. The authors studied in \cite{bbg3} the interaction of a biological species and a chemical substance of non-diffusive nature and a chemotactic model with non-local terms in \cite{bbg4}.

The paper is organized as follows: in Section 2 we explain the derivation of the explicit formulae method. In Section 3 we study the convergence of the explicit scheme. Next, in Section 4 we give several examples over regular and irregular meshes. Some conclusions are finally drawn.
\section{Preliminaries of de GFD Method}\label{preliminares}
\begin{figure}[htp]
	\centering
	\includegraphics[scale=0.45]{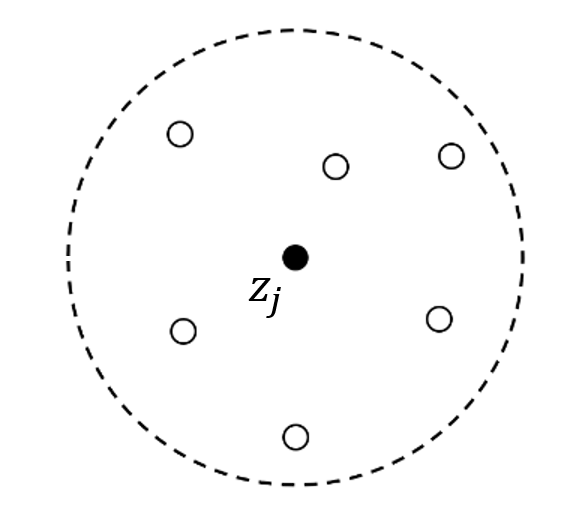}
	\caption{Representation of the nodes that make up an $E_6-$star centered in $\vect{z_j}$.}
	\label{fig:4.1}
\end{figure}

We first introduce the preliminaries of the Generalized Finite Difference Method. For this purpose, let us consider the following general problem over a bounded domain, $\Omega \subset \mathbb{R}^N$:
\begin{equation} \label{4.1}
\begin{cases}
\displaystyle  \frac{\partial u}{\partial t} (\vect{x},t) = L_\Omega [u(\vect{x},t)], \quad \text{for } \vect{x} \in \Omega, \hspace{0.1 cm} t>0, \\\\
\displaystyle  \frac{\partial u}{\partial \vec{n}}  (\vect{x},t) = 0, \quad \text{for } \vect{x} \in \partial \Omega, \hspace{0.1 cm}  t>0, \\\\
\displaystyle  u(\vect{x},0) = g(\vect{x}), \quad \text{for } \vect{x} \in \Omega,
\end{cases}
\end{equation}
where $L_\Omega$ is a nonlinear differential operator defined in $\Omega$. In all the subsequent development, for the sake of simplicity in the graphical and analytical representation, we take $N = 2$, but higher dimensions can be analogously treated.
\\\\
Hence, we consider a discretization of $\Omega$, consisting of a set of points which we denote by $M = \{z_1, \dots, z_m\} \subset \Omega$, known as \textit{nodes}, each of them being of the form $\vect{z_j} = (x_j,y_j)$.  As in the standard Finite Difference method, the objective here is to obtain an approximation of the solution of (\ref{4.1}) at each of the nodes of $M$, using in the same way finite difference formulas evaluated on them. The difference between the two methods, however, lies on the fact that for Generalized Finite Differences the nodes do not need to form a regular mesh over $\Omega$ and can be randomly distributed instead, though some considerations shall be made. The generalized finite
difference method was created in the 1970s with the work of Jensen \cite{jensen},
Perrone and Kao \cite{perrone}, although it was already suggested by Collatz \cite{collatz} and
Forsythe and Wasow \cite{forsythe} ten years earlier, and a large number of authors
have contributed to its progress and improvement, as can be seen in \cite{albuquerque}.
\\\\
Moreover, with \eqref{4.1} being an evolution problem, a temporal discretization must also be taken into account, for which we take a constant time step, $\Delta t>0$. As per usual, let $U_j^n$ be the approximation of $u(\vect{z_j}, n \Delta t)$, whose value we seek to determine. Using a progressive difference formula, we consider the first order approximation given by:
\begin{equation} \label{4.2}
\frac{\partial u}{\partial t}(\vect{z_j}, n \Delta t) \approx \frac{U_j^{n+1} - U_j^n}{\Delta t}.
\end{equation}
For the spatial derivatives, we assume $L_\Omega$ to be a second order operator, that is, only including $\frac{\partial }{\partial x}$, $\frac{\partial }{\partial y}$, $\frac{\partial^2 }{\partial x^2}$, $\frac{\partial^2 }{\partial y^2}$ y $\frac{\partial^2 }{\partial x \partial y}$. For each $\vect{z_j} \in M$, we consider a set of $s$ other nodes, located within a neighbourhood of $\vect{z_j}$, called an $E_s -$star, as shown in Figure \ref{fig:4.1}.
\\\\On this $E_s-$star centered in $\vect{z_j} = (x_j,y_j)$ we define the following distances for the other nodes $\vect{z_i} = (x_i, y_i) \in E_s$:
$$h_i^j = x_i - x_j, \quad k_i^j = y_i - y_j,$$
and through a second order Taylor approximation, we have:
\begin{equation} \label{4.2.1}
\begin{split}
&u(x_i, y_i) \approx u(x_j, y_j) + h_i^j \frac{\partial u}{\partial x}(x_j, y_j) + k_i^j \frac{\partial u}{\partial y}(x_j, y_j) \\&+ \frac{1}{2}\left[(h_i^j)^2 \frac{\partial^2 u}{\partial x^2}(x_j, y_j) +  (k_i^j)^2 \frac{\partial^2 u}{\partial y^2}(x_j, y_j) + 2h_i^j k_i^j \frac{\partial^2 u}{\partial x \partial y}(x_j, y_j)  \right].
\end{split}
\end{equation} 
Thus, the following function $B$ can be considered, representing a weighted sum of quadratic errors between the value of $U_i^n$ and its second order Taylor approximation \eqref{4.2.1} centered in $\vect{z_j}$:
\begin{equation} \label{4.3}
\begin{split}
B(U_j^n) =& \sum_{i=1}^s (w_i^j)^2 \Biggl[ U_j^n + h_i^j \frac{\partial U_j^n}{\partial x} +  k_i^j \frac{\partial U_j^n}{\partial y} \\&+ \frac{1}{2} \left( (h_i^j)^2 \frac{\partial^2  U_j^n}{\partial x^2} +  (k_i^j)^2 \frac{\partial^2  U_j^n}{\partial y^2} + 2h_i^j k_i^j \frac{\partial^2  U_j^n}{\partial x \partial y} \right) - U_i^n  \Biggr]^2.
\end{split}
\end{equation}
Here, the partial derivatives of $U_j^n$ are a symbolic notation referring to their approximations, which we look to determine. The weights $w_i^j$ are non-negative symmetric functions that decrease with the distance from $\vect{z_j}$ to $\vect{z_i}$, usually given by:
$$w_i^j = \frac{1}{(h_i^j + k_i^j)^{\alpha/2}} = \frac{1}{||\vect{z_j} - \vect{z_i}||^\alpha},$$
for a certain $\alpha > 0$. As a result, to find the best approximations formulas for the partial derivatives, $B$ has to be minimized with respect to $\partial U_j^n/\partial x$, $\partial U_j^n/\partial y$, $\partial^2 U_j^n/\partial x^2$, $\partial^2 U_j^n/\partial y^2$ and $\partial^2 U_j^n/\partial x \partial y$. Therefore, computing the partial derivatives of $B$ with respect to these variables and setting them to zero, we arrive at a linear system of equations, that can be expressed as:
\begin{equation} \label{4.4}
\vect{A_j^n} \cdot \vect{D_j^n} = \vect{b_j^n},
\end{equation}
where:
\begin{equation} \label{4.5}
\vect{A_j^n} = \left(
\begin{array}{cccc}
h^j_1 & h^j_2 & \cdots & h^j_s \\
k^j_1 & k^j_2 & \cdots & k^j_s \\
\vdots & \vdots & \vdots & \vdots \\
h^j_1k^j_1 & h^j_2k^j_2 & \cdots & h^j_sk^j_s \\
\end{array}
\right)\left(
\begin{array}{cccc}
(w_1^j)^2 &  &  &  \\
& (w_2^j)^2 &  &  \\
&  & \cdots &  \\
&  &  & (w_s^j)^2 \\
\end{array}
\right)\left(
\begin{array}{cccc}
h^j_1 & k^j_1 & \cdots & h^j_1k^j_1 \\
h^j_2 & k^j_2 & \cdots & h^j_2k^j_2 \\
\vdots & \vdots& \vdots & \vdots \\
h^j_s & k^j_s &  \cdots & h^j_sk^j_s \\
\end{array}\right)  
%
%
%
%
%
%
\end{equation}
\vspace{0.02 cm}
\begin{equation} \label{4.6}
\vect{b_j^n} = \begin{pmatrix}
\displaystyle    \sum_{i=1}^s (U_i^n - U_j^n)h_i^j (w_i^j)^2 \\\\ \displaystyle \sum_{i=1}^s (U_i^n - U_j^n)k_i^j (w_i^j)^2 \\\\ \displaystyle \frac{1}{2} \sum_{i=1}^s (U_i^n - U_j^n)(h_i^j)^2 (w_i^j)^2 \\\\ \displaystyle
\frac{1}{2} \sum_{i=1}^s (U_i^n - U_j^n)(k_i^j)^2 (w_i^j)^2 \\\\ \displaystyle
\sum_{i=1}^s (U_i^n - U_j^n)h_i^j k_i^j (w_i^j)^2
\end{pmatrix},
\end{equation}
and $\vect{D_j^n}$ is the unknown vector, containing the partial derivatives:
\begin{equation} \label{4.7}
\vect{D_j^n} = \begin{pmatrix} \displaystyle 
\frac{\partial U_j^n}{\partial x} & \displaystyle \frac{\partial U_j^n}{\partial y} &  \displaystyle \frac{\partial^2 U_j^n}{\partial x^2} &  \displaystyle \frac{\partial^2 U_j^n}{\partial y^2} &  \displaystyle \frac{\partial^2 U_j^n}{\partial x \partial y}
\end{pmatrix}^T.
\end{equation}
Thus, the solution to system \eqref{4.4} yields the desired approximation formulas for the partial derivatives, based only on the approximate values of $u$ on the nodes of the star. To avoid possible confusions caused by the super and sub indexes, without loss of generality we center the study on a certain node, denoted by $\vect{z_0}$, taking into account the following notation simplification:
\begin{itemize}
	\item The node $\vect{z_j} = (x_j, y_j)$ is substituted by $\vect{z_0} = (x_0, y_0)$.
	\item The matrix $\vect{A_j^n}$ is only denoted by $\vect{A}$.
	\item The distances $h_i^j$ y $k_i^j$ and weights $w_i^j$ simply become $h_i$, $k_i$ y $w_i$.
\end{itemize}
Moreover, let us define for each node of the $E_s-$star centered in $\vect{z_0}$ the vector:
\begin{equation} \label{4.7.2}
\vect{c_i} = \begin{pmatrix}
h_i & k_i & \frac{1}{2} h_i^2 &\frac{1}{2} k_i^2 &  h_i k_i
\end{pmatrix},
\end{equation}
such that the solution to system \eqref{4.4} can be expressed as:
\begin{equation} \label{4.8}
\begin{cases}
\displaystyle   \frac{\partial U(x_0,y_0, n \Delta t)}{\partial x} = - \lambda_{01} U_0^n + \sum_{i=1}^s \lambda_{i1}U_i^n + \mathcal{O}(h_i^2,k_i^2), \\\\\
\displaystyle \frac{\partial U(x_0,y_0, n \Delta t)}{\partial y} = - \lambda_{02} U_0^n + \sum_{i=1}^s \lambda_{i2}U_i^n +  \mathcal{O}(h_i^2,k_i^2), \\\\
\displaystyle \frac{\partial^2 U(x_0,y_0, n \Delta t)}{\partial x^2} = - \lambda_{03} U_0^n + \sum_{i=1}^s \lambda_{i3}U_i^n +  \mathcal{O}(h_i^2,k_i^2), \\\\
\displaystyle \frac{\partial^2 U(x_0,y_0, n \Delta t)}{\partial y^2} = - \lambda_{04} U_0^n + \sum_{i=1}^s \lambda_{i4}U_i^n +  \mathcal{O}(h_i^2,k_i^2), \\\\
\displaystyle \frac{\partial^2 U(x_0,y_0, n \Delta t)}{\partial x \partial y} = - \lambda_{05} U_0^n + \sum_{i=1}^s \lambda_{i5}U_i^n +  \mathcal{O}(h_i^2,k_i^2),
\end{cases}
\end{equation}
where the coefficients $\lambda_{i,j}$ are given by:
\begin{equation} \label{4.9}
\lambda_{ir} = w_i^2 (\vect{A}^{-1} \vect{c_i})_r, \quad \lambda_{0i} = \sum_{i=1}^s \lambda_{ir},
\end{equation}
for $r\in \{1, \dots, s\}$, $i \in \{1, \dots, 5\}$, being $(\vect{A}^{-1} \vect{c_i})_r$ the $r-$th coordinate of the vector $\vect{A}^{-1} \vect{c_i}$. Hence, combining the finite difference formulas given by equations \eqref{4.8} with the forward differences approximation for the temporal derivative from \eqref{4.2}, an explicit scheme for solving \eqref{4.1} is obtained. As a remark, due to the fact that $\vect{A}$ is a positive definite matrix (see, for example \cite{R6,vargas}), $\vect{A}^{-1}$ can be computed through a Cholesky factorization.
\section{Convergence of the numerical method}\label{metodo}
In this section, we present a numerical scheme for solving problem \eqref{1.2} through the Generalized Finite Difference Method, proving its convergence. To do so, it suffices to use the approximation formulae obtained in Section \ref{preliminares} for the partial derivatives, taking into account that \eqref{1.2} is a system of two equations. Therefore, substituting the derivatives by the finite difference formulae and rearranging the terms, for every inner node we obtain:
\begin{equation}
\label{4.13}
\left\lbrace
\begin{aligned}{}
&\frac{U_0^{n+1}-U_0^n}{\Delta t}=\gamma(V^{n}_0)\left[-\lambda_{00}U_0^n+\sum^{s}_{i=1}\lambda_{i0}U_i^n\right] \\&\qquad\qquad
+2\gamma'(V^n_0)\left(-\lambda_{01}U_0^n+\sum^{s}_{i=1}\lambda_{i1}U_i^n\right)\left(-\lambda_{01}V_0^n+\sum^{s}_{i=1}\lambda_{i1}V_i^n\right)
\\&\qquad\qquad  +2\gamma'(V^n_0)\left(-\lambda_{02}U_0^n+\sum^{s}_{i=1}\lambda_{i2}U_i^n\right)\left(-\lambda_{02}V_0^n+\sum^{s}_{i=1}\lambda_{i2}V_i^n\right) 
\\&\qquad\qquad+ U^n_0\gamma''(V^n_0)\left[\left(-\lambda_{01}V_0^n+\sum^{s}_{i=1}\lambda_{i1}V_i^n\right)^2+\left(-\lambda_{02}V_0^n+\sum^{s}_{i=1}\lambda_{i2}V_i^n\right)^2\hspace{0.05 cm}\right]\\
&\qquad\qquad+ U^n_0\gamma'(V^n_0)\Big( V^n_0 -U^n_0 \Big) + \mu U^n_0(1-U^n_0) + \mathcal{O}(\Delta t,h_i^2,k_i^2),
\end{aligned}
\right.
\end{equation} and
\begin{equation}\label{4.14}
V_0^n-\left[-\lambda_{00}V_0^n+\sum^{s}_{i=1}\lambda_{i0}V_i^n\right]=U^n_0 +\mathcal{O}(\Delta t,h_i^2,k_i^2).
\end{equation}
Notice that since the second equation of \eqref{1.2} is elliptic, its subsequent scheme is implicit, as no time derivative is present. Algorithmically, the method is initialized with the values of $u_0(\vect{x})$, then obtaining $V_0^0$ for all inner nodes through the solution to the linear system resulting from \eqref{4.14}. The value of $U_0^1$ can then be directly computed through \eqref{4.13}, consequently calculating $V_0^1$ and repeating this process until the desired final time.
\\\\The main result of the section is the proof of the convergence of the scheme \eqref{4.13}-\eqref{4.14}, stated in the following theorem.
\begin{thm} \label{num1}
	Let $(u,v)$ be the solution to system \eqref{1.2} in the sense specified in \cite{tello} under conditions \eqref{1.4}-\eqref{1.5}, then the GFD scheme \eqref{4.13}-\eqref{4.14} is convergent if the time step $\Delta t$ is taken such that the following bound holds for every inner node:
	$$\Delta t<\dfrac{2+|\lambda_{00}|+|\sum_{i=1}^s|\lambda_{i0}|}{\Biggl[|1-\lambda_{00}|+\sum_{i=1}^s|\lambda_{i0}|\Biggr](A'_1+A''_1)+B_1},$$
	where the explicit expressions of the coefficients $A'_1$, $A''_1$ and $B_1$ are given in equations in the proof.
\end{thm}
Proof:
	Since the exact values $u,v$ must fulfil system (\ref{1.2}), we take the difference between the approximated solution given by the GFD scheme and the exact expression of the continuous solution. For simplicity, let us denote by $eu^n_0:=U^n_0-u^n_0$ the difference of the discrete and the continuous solution at the node $\vect{z_0}$ (and in the same manner we define $eu^n_i,$ and $ev^n_i$). By the great symmetry of the scheme, let us only perform the computations explicitly for the most significant terms, as the rest are analogously treated.
	\\\\Firstly, by the Mean Value Theorem it yields:
	\begin{equation}\label{prueba}
	\begin{split}
	\gamma&(V^n_0)\left[-\lambda_{00}U_0^n+\sum^{s}_{i=1}\lambda_{i0}U_i^n\right]-\gamma(v^n_0)\left[-\lambda_{00}u_0^n+\sum^{s}_{i=1}\lambda_{i0}u_i^n\right] \\&\pm\gamma(V^n_0)\left[-\lambda_{00}u_0^n+\sum^{s}_{i=1}\lambda_{i0}u_i^n\right]=\\
	&=\gamma'(\xi)\left[-\lambda_{00}u_0^n+\sum^{s}_{i=1}\lambda_{i0}u_i^n\right]ev^n_0+\gamma(V^n_0)\left[-\lambda_{00}eu_0^n+\sum^{s}_{i=1}\lambda_{i0}eu_i^n\right].
	\end{split}
	\end{equation}
	Then, in a similar way, the following holds:
	\begin{equation}\label{prueba1}
	\begin{split}
	\small &2\gamma'(V^n_0)\left(-\lambda_{01}U_0^n+\sum^{s}_{i=1}\lambda_{i1}U_i^n\right)\left(-\lambda_{01}V_0^n+\sum^{s}_{i=1}\lambda_{i1}V_i^n\right)\\&-2\gamma'(v^n_0)\left(-\lambda_{01}u_0^n+\sum^{s}_{i=1}\lambda_{i1}u_i^n\right)\left(-\lambda_{01}v_0^n+\sum^{s}_{i=1}\lambda_{i1}v_i^n\right)\\
	&\pm 2\gamma'(V^n_0)\left(-\lambda_{01}u_0^n+\sum^{s}_{i=1}\lambda_{i1}u_i^n\right)\left(-\lambda_{01}v_0^n+\sum^{s}_{i=1}\lambda_{i1}v_i^n\right)\\& \pm 2\gamma'(V^n_0)\left(-\lambda_{01}u_0^n+\sum^{s}_{i=1}\lambda_{i1}u_i^n\right)\left(-\lambda_{01}V_0^n+\sum^{s}_{i=1}\lambda_{i1}V_i^n\right)=\\
	&=2\gamma''(\xi)\left(-\lambda_{01}u_0^n+\sum^{s}_{i=1}\lambda_{i1}u_i^n\right)\left(-\lambda_{01}v_0^n+\sum^{s}_{i=1}\lambda_{i1}v_i^n\right)ev^n_0 \\&+2\gamma'(V^n_0)\left(-\lambda_{01}u_0^n+\sum^{s}_{i=1}\lambda_{i1}u_i^n\right)\left(-\lambda_{01}ev_0^n+\sum^{s}_{i=1}\lambda_{i1}ev_i^n\right)+\\
	&2\gamma'(V^n_0)\left(-\lambda_{01}eu_0^n+\sum^{s}_{i=1}\lambda_{i1}eu_i^n\right)\left(-\lambda_{01}V_0^n+\sum^{s}_{i=1}\lambda_{i1}V_i^n\right).
	\end{split}
	\end{equation}
	After applying the same treatment to the rest of the terms, we obtain:
	\begin{equation}\label{prueba7}
	\begin{split}
	\small &\frac{eu^{n+1}_0-eu^n_0}{\Delta t}=\gamma''(\xi)\left[-\lambda_{00}u_0^n+\sum^{s}_{i=1}\lambda_{i0}u_i^n\right]ev^n_0+\gamma'(V^n_0)\left[-\lambda_{00}eu_0^n+\sum^{s}_{i=1}\lambda_{i0}eu_i^n\right]\\
	&+2\gamma''(\xi)\left(-\lambda_{01}U_0^n+\sum^{s}_{i=1}\lambda_{i1}U_i^n\right)\left(-\lambda_{01}V_0^n+\sum^{s}_{i=1}\lambda_{i1}V_i^n\right)ev^n_0+\\
	&+2\gamma'(V^n_0)\left(-\lambda_{01}u_0^n+\sum^{s}_{i=1}\lambda_{i1}u_i^n\right)\left(-\lambda_{01}ev_0^n+\sum^{s}_{i=1}\lambda_{i1}ev_i^n\right)+\\
	&+2\gamma'(V^n_0)\left(-\lambda_{01}eu_0^n+\sum^{s}_{i=1}\lambda_{i1}eu_i^n\right)\left(-\lambda_{01}V_0^n+\sum^{s}_{i=1}\lambda_{i1}V_i^n\right)+\\
	&+2\gamma''(\xi)\left(-\lambda_{02}U_0^n+\sum^{s}_{i=1}\lambda_{i2}U_i^n\right)\left(-\lambda_{02}V_0^n+\sum^{s}_{i=1}\lambda_{i2}V_i^n\right)ev^n_0+\\
	&+2\gamma'(V^n_0)\left(-\lambda_{02}eu_0^n+\sum^{s}_{i=1}\lambda_{i2}eu_i^n\right)\left(-\lambda_{02}V_0^n+\sum^{s}_{i=1}\lambda_{i2}V_i^n\right)+\\
	&+2\gamma'(V^n_0)\left(-\lambda_{02}eu_0^n+\sum^{s}_{i=1}\lambda_{i2}eu_i^n\right)\left(-\lambda_{02}V_0^n+\sum^{s}_{i=1}\lambda_{i2}V_i^n\right)+\\
	&+eu^n_0\gamma''(V^n_0)\Biggl(-\lambda_{01}V_0^n+\sum^{s}_{i=1}\lambda_{i1}V_i^n\Biggr)^2+eu^n_0\gamma''(V^n_0)\Biggl(-\lambda_{02}V_0^n+\sum^{s}_{i=1}\lambda_{i2}V_i^n\Biggr)^2\\
	&u^n_0\gamma'''(\xi)\Biggl(-\lambda_{01}V_0^n+\sum^{s}_{i=1}\lambda_{i1}V_i^n\Biggr)^2ev^n_0+u^n_0\gamma'''(\xi)\Biggl(-\lambda_{02}V_0^n+\sum^{s}_{i=1}\lambda_{i2}V_i^n\Biggr)^2ev^n_0\\
	&+u^n_0\gamma''(V^n_0)\Biggl(-\lambda_{01}ev_0^n+\sum^{s}_{i=1}\lambda_{i1}ev_i^n\Biggr)\Biggl(-\lambda_{01}(v_0^n+V^n_0)+\sum^{s}_{i=1}\lambda_{i1}(v_i^n+V^n_i)\Biggr)\\
	&+u^n_0\gamma''(V^n_0)\Biggl(-\lambda_{02}ev_0^n+\sum^{s}_{i=1}\lambda_{i2}ev_i^n\Biggr)\Biggl(-\lambda_{02}(v_0^n+V^n_0)+\sum^{s}_{i=1}\lambda_{i2}(v_i^n+V^n_i)\Biggr)\\
	&+eu^n_0V^n_0\gamma'(V^n_0)+u^n_0v^n_0\gamma''(\xi)ev^n_0+u^n_0ev^n_0\gamma'(V^n_0)-(U^n_0)^2\gamma''(\xi)ev^n_0\\
	&-eu^n_0(u^n_0+U^n_0)\gamma'(v^n_0)+\mu eu^n_0-\mu eu^n_0(u^n_0+U^n_0)
	\end{split}
	\end{equation}
	Denoting $eu^n=\displaystyle\max_{i=0,...,s}\{|eu^n_i|\}$ y $ev^n = \displaystyle\max_{i=0,...,s}\{|ev^n_i|\}$ and taking bounds in the last expression, we get:
	\begin{equation}\label{prueba8}
	eu^{n+1}\leq A_1eu^n+B_1ev^n,
	\end{equation}
	where the positive coefficients $A_1,$ y $ B_1$ are given by:
	\begin{equation}\label{prueba9}
	\begin{split}
	A_1:&=\Biggl|1+ \Delta t\Biggl[-\lambda_{00}-2\gamma'(V^n_0)\lambda_{01}\Biggl(-\lambda_{01}V^n_0+\sum_{i=1}^s\lambda_{i1}V^n_i\Biggr)\\
	&-2\gamma'(V^n_0)\lambda_{02}\Biggl(-\lambda_{02}V^n_0+\sum_{i=1}^s\lambda_{i2}V^n_i\Biggr)+\gamma''(V^n_0)\Biggl(-\lambda_{01}V^n_0+\sum_{i=1}^s\lambda_{i1}V^n_i\Biggr)^2\\
	&+\gamma''(V^n_0)\Biggl(-\lambda_{02}V^n_0+\sum_{i=1}^s\lambda_{i2}V^n_i\Biggr)^2+V^n_0\gamma'(V^n_0)+u^n_0\gamma'(V^n_0)\\
	&-(u^n_0+U^n_0)\gamma'(v^n_0)+\mu-\mu(u^n_0+U^n_0)\Biggr]\Biggr|+\Delta t\Biggl[|\gamma'(v^n_0)\sum_{i=1}^s\lambda_{i0}|\\
	&+2|\gamma''(v^n_0)\Biggl(-\lambda_{01}V^n_0+\sum_{i=1}^s\lambda_{i1}V^n_i\Biggr)\sum_{i=1}^s\lambda_{i1}|\\&+2|\gamma''(v^n_0)\Biggl(-\lambda_{02}V^n_0+\sum_{i=1}^s\lambda_{i2}V^n_i\Biggr)\sum_{i=1}^s\lambda_{i2}|\Biggr],
	\end{split}
	\end{equation}
	\vspace{0.3 cm}
	\begin{equation}\label{prueba10}
	\begin{split}
	B_1:&=\Delta t\Biggl|\gamma'(\xi)\Biggl(-\lambda_{00}U^n_0+\sum_{i=1}^s\lambda_{i0}U^n_i\Biggr)\\&+2\gamma''(\xi)\Biggl(-\lambda_{01}U^n_0+\sum_{i=1}^s\lambda_{i1}U^n_i\Biggr)\Biggl(-\lambda_{01}V^n_0+\sum_{i=1}^s\lambda_{i1}V^n_i\Biggr)\\
	&+2\gamma''(\xi)\Biggl(-\lambda_{02}U^n_0+\sum_{i=1}^s\lambda_{i2}U^n_i\Biggr)\Biggl(-\lambda_{02}V^n_0+\sum_{i=1}^s\lambda_{i2}V^n_i\Biggr)\\&-2\gamma'(V^n_0)\Biggl(-\lambda_{01}U^n_0+\sum_{i=1}^s\lambda_{i1}U^n_i\Biggr)\lambda_{01}
	-2\gamma'(V^n_0)\Biggl(-\lambda_{02}U^n_0+\sum_{i=1}^s\lambda_{i2}U^n_i\Biggr)\lambda_{02}\\&+u^n_0\gamma'''(\xi)\Biggl(-\lambda_{01}V^n_0+\sum_{i=1}^s\lambda_{i1}V^n_i\Biggr)^2 +u^n_0\gamma'''(\xi)\Biggl(-\lambda_{02}V^n_0+\sum_{i=1}^s\lambda_{i2}V^n_i\Biggr)^2 \\&-u^n_0\gamma''(v^n_0)\Biggl(-\lambda_{01}(v^n_0+V^n_0)+\sum_{i=1}^s\lambda_{i1}(v^n_i+V^n_i)\Biggr)\lambda_{01}\\
	&-u^n_0\gamma''(v^n_0)\Biggl(-\lambda_{02}(v^n(_0+V^n_0)+\sum_{i=1}^s\lambda_{i2}(v^n_i+V^n_i)\Biggr)\lambda_{02}\\&+u^n_0v^n_0\gamma''(\xi)-(U^n_0)^2\gamma''(\xi)\Biggr|+\Delta t\Biggl[2|\gamma'(v^n_0)\Biggl(-\lambda_{01}U^n_0+\sum_{i=1}^s\lambda_{i1}U^n_i\Biggr)\sum_{i=1}^s\lambda_{i1}|\\&+2|\gamma'(v^n_0)\Biggl(-\lambda_{02}U^n_0+\sum_{i=1}^s\lambda_{i2}U^n_i\Biggr)\sum_{i=1}^s\lambda_{i2}|\\&+|u^n_0\gamma''(v^n_0)\Biggl(-\lambda_{01}(v^n(_0+V^n_0)+\sum_{i=1}^s\lambda_{i1}(v^n_i+V^n_i)\Biggr)\sum_{i=1}^s\lambda_{i1}|\\&+|u^n_0\gamma''(v^n_0)\Biggl(-\lambda_{02}(v^n(_0+V^n_0)+\sum_{i=1}^s\lambda_{i2}(v^n_i+V^n_i)\Biggr)\sum_{i=1}^s\lambda_{i2}| \Biggr].
	\end{split}
	\end{equation}
	We can now rewrite $A_1$ as $A_1=|1-\Delta t A_1^{'}|+\Delta t A_{1}^{''},$ for an obvious choice of $A_1^{'}$ and $A_1^{''}$.
	\\\\Following a similar procedure with \eqref{4.14}, subtracting the exact solution and taking bounds, we have: 
	\begin{equation}\label{prueba13}
	eu^{n} \leq \Biggl[|1-\lambda_{00}|+\sum_{i=1}^s|\lambda_{i0}|\Biggr] ev^n.
	\end{equation}
	By substitution of \eqref{prueba13} in \eqref{prueba8} for $eu^{n}$ we arrive to
	\begin{equation}\label{prueba13b}
	eu^{n+1}\leq \Biggl(A_1\Biggl[|1-\lambda_{00}|+\sum_{i=1}^s|\lambda_{i0}|\Biggr]+B_1\Biggr)ev^n.
	\end{equation}
	 Since we aim to obtain a convergent scheme, for each star we have to assure that $eu^{n}$ and $ev^{n}$ tend to 0. We impose
	\begin{equation} \label{prueba14}
	A_1[|1-\lambda_{00}|+\sum_{i=1}^s|\lambda_{i0}|]+B_1 < 1. 
	\end{equation}
	Having in mind $A_1=|1-\Delta t A_1^{'}|+\Delta t A_{1}^{''}$ $A_1'$, we want to establish an upper bound for  $\Delta t$ from \eqref{prueba14}, so it follows:
	\begin{equation} \label{prueba15}
	\Delta t<\dfrac{2+|\lambda_{00}|+|\sum_{i=1}^s|\lambda_{i0}|}{\Biggl[|1-\lambda_{00}|+\sum_{i=1}^s|\lambda_{i0}|\Biggr](A'_1+A''_1)+B_1}.
	\end{equation}
\section{Numerical Tests}
Lastly, in this section we use the GFD scheme \eqref{4.13}-\eqref{4.14} to numerically solve system \eqref{1.2}, in order to illustrate the asymptotic behavior of its solution. 
\subsection{Example 1}
\begin{figure}[htp]
	\centering
	\includegraphics[scale=0.5]{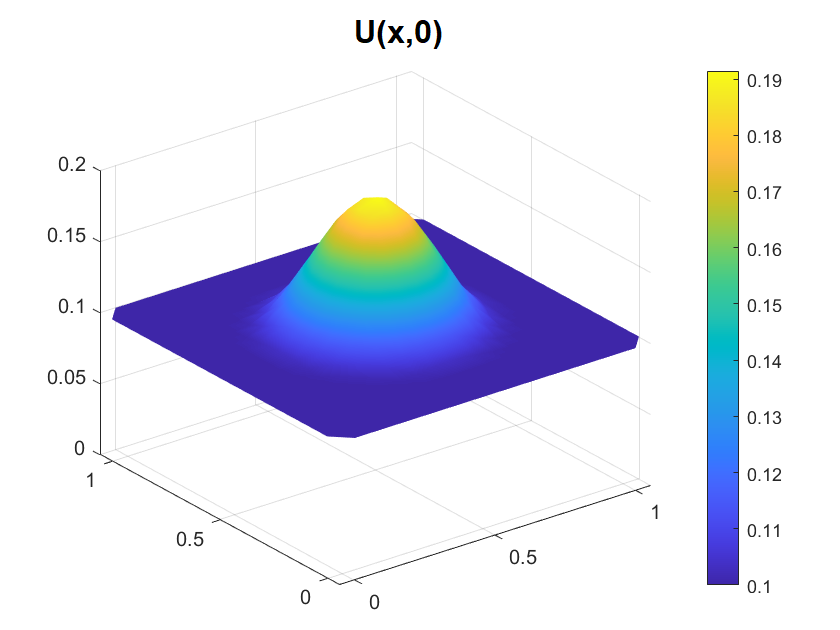}
	\caption{Initial bacterial density, $u(\vect{x},0)$ from Example 1.}
	\label{fig:4.1.2}
\end{figure}

\begin{figure}[htp]
	\centering
	\includegraphics[scale=0.3]{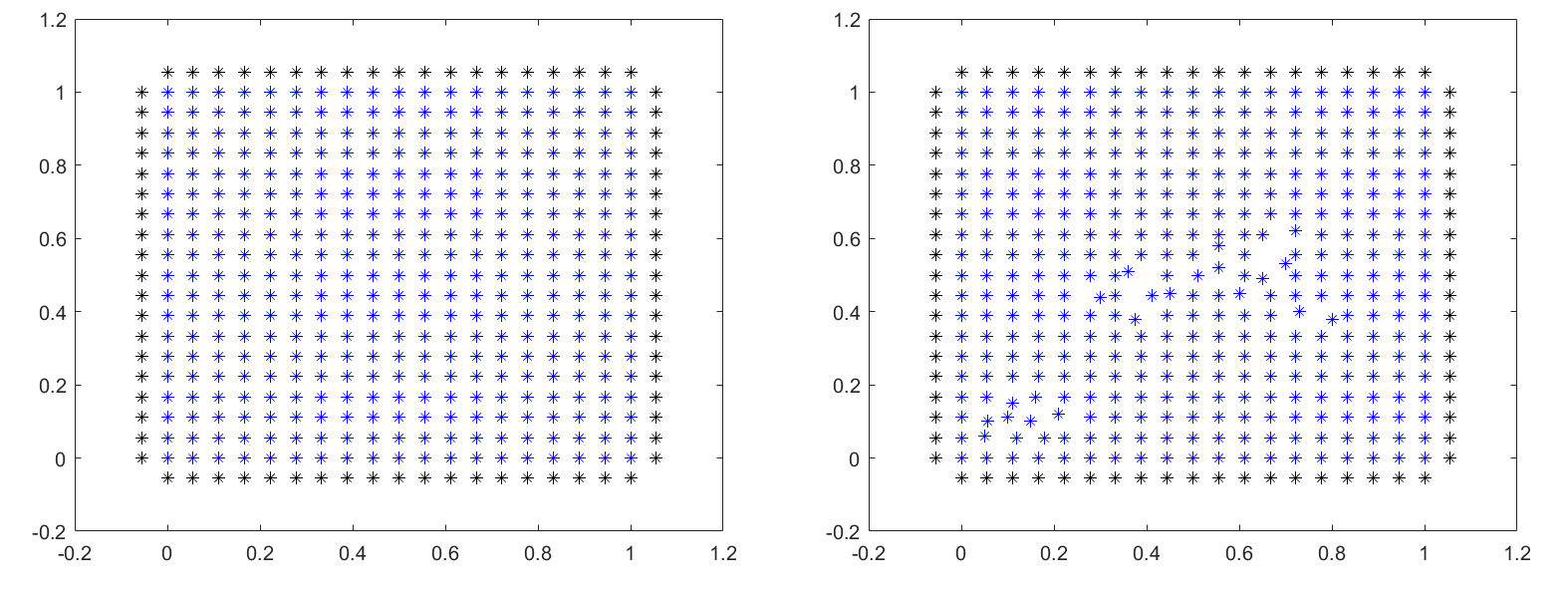}
	\caption{Discretizations of $\Omega$ considered.}
	\label{fig:4.2}
\end{figure}

\begin{figure}[htp]
	\centering
	\includegraphics[scale=0.3]{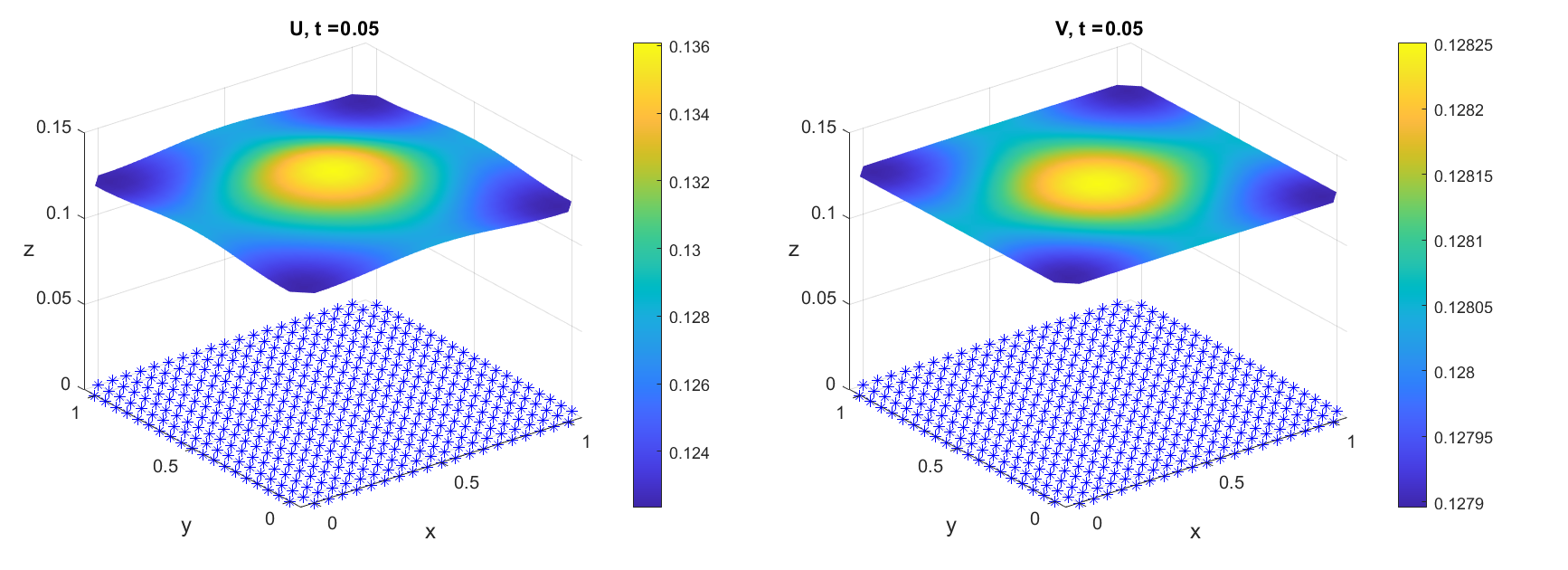}
	\caption{Numerical solution $(U,V)$ from Example 1 at $t = 0.05$ using the regular grid.}
	\label{fig:4.3}
\end{figure}

\begin{figure}[htp]
	\centering
	\includegraphics[scale=0.45]{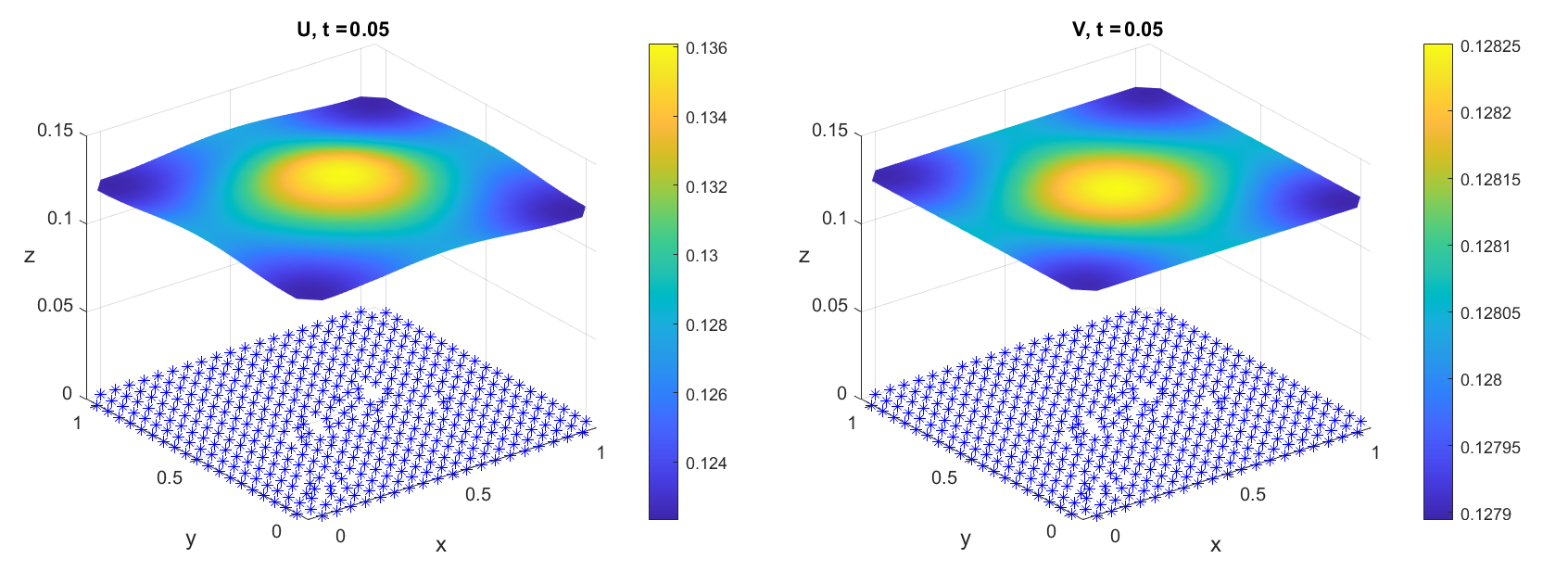}
	\caption{Numerical solution $(U,V)$ at $t = 0.05$ using the irregular cloud of points.}
	\label{fig:4.4}
\end{figure}

We take the square domain $\Omega = [0,1]\times [0,1]$ and as initial bacterial density, we consider the following radially-symmetric function:
\begin{equation}\label{5.1}
      u(\vect{x},0) = a + b \hspace{0.05 cm}\varphi(r), \hspace{0.2 cm} \text{with} \hspace{0.2 cm} \varphi(r) = \begin{cases}
      \exp \left(\frac{-1}{0.25 - r^2} \right), \hspace{0.2 cm} \text{if} \hspace{0.2 cm} r < 0.5, \\
      0,\hspace{0.2 cm} \text{otherwise},
    \end{cases},
\end{equation}
where $r := (x-1/2)^2 + (y-1/2)^2$, the euclidean distance to the center of $\Omega$. In \eqref{5.1}, $a$, $b > 0$, are positive constants that provide shape and scale to the function. We set $a = 0.1$ and $b = 5$, and thus this choice of $u(\vect{x},0)$ satisfies hypothesis \eqref{1.5}. The graph of $u(\vect{x},0)$ is shown in Figure \ref{fig:4.1.2}
For the numerical method, we consider two different discretizations of $\Omega$, the first one consisting of a regular grid of $19 \times 19$ points, while the second one is the irregular cloud of points depicted in Figure \ref{fig:4.2}.
The surrounding black nodes in both discretizations are fictitious nodes, employed for approximating the normal derivatives $\frac{\partial u}{\partial \Vec{n}}$ and $\frac{\partial v}{\partial \Vec{n}}$ in $\partial \Omega$ through a first order scheme for the Neumann homogeneous conditions. 
\\\\
Regarding $\gamma$, the motility regulation function, we take $\gamma(v) := e^{-v}$, that fulfils \eqref{1.4} for a large enough $\mu$, the growth factor of the logistic model. In particular, as we need to verify $-2 \gamma'(s) + \gamma''(s) s  < \mu$ for all $s \geq 0$, for this choice of $\gamma$ we have $-2 \gamma'(s) + \gamma''(s) s = 2e^{-s} + se^{-s}$, which has a maximum at $s = -1$ with a value of $e$. Hence, we can consider $\mu = 3$. We set the weights to be $w_{i} := \frac{1}{h_i^2 + k_i^2}$ and $\Delta t = 0.001$, satisfying the assumption made in Theorem \ref{num1}. 
\\\\
Figures \ref{fig:4.3} and \ref{fig:4.4} show the numerical solution $U$, $V$ of system \eqref{1.2} on $t = 0.05$, obtained by implementing the GFD scheme \eqref{4.13}-\eqref{4.14} through a MATLAB code. The values depicted in Figure \ref{fig:4.3} are calculated using the regular grid, whereas Figure \ref{fig:4.4} shows the results for the irregular cloud of points.

Moreover, we present the values of $||U-1||_{l^\infty(\Omega)}$ and $||V-1||_{l^\infty(\Omega)}$ in various time instants in Table \ref{t4.1}, calculated using the regular grid. As we can see, this verifies the convergence \eqref{1.4.1}
$$\lim_{t\to\infty} ||u-1||_{L^\infty(\Omega)} + ||v-1||_{L^\infty(\Omega)} = 0,$$
with the bacterial species growing to their maximum carrying capacity and the AHL reaching the homogeneous state $\alpha$, obtained after scaling the variables.
\begin{table}[h]
\centering
\begin{tabular}{llllll}\hline \\[-1em]
$T(s)$  & 0.05 & 1 & 2.5 & 5 & 10 \\ \\[-1em] \hline
\\[-0.5em]
$||U-1||_{l^\infty(\Omega)}$   & 0.8777 & 0.2821 & 0.0043  & $2.3740 \cdot 10^{-6}$ & $6.3094 \cdot 10^{-13}$  \\ \\[-0.5em]
$||V-1||_{l^\infty(\Omega)}$    & 0.8721 & 0.2827 & 0.0043  & $2.3811 \cdot 10^{-6}$ & $2.3438 \cdot 10^{-12}$  \\ \\[-0.8em] \hline
\end{tabular}
\caption{\label{t4.1} Values of $||U-1||_{l^\infty(\Omega)}$ and $||V-1||_{l^\infty(\Omega)}$ in different time instants from Example 1, calculated with the regular grid.}
\end{table}
\subsection{Example 2}
\begin{figure}[htp]
    \centering
    \includegraphics[scale=0.3]{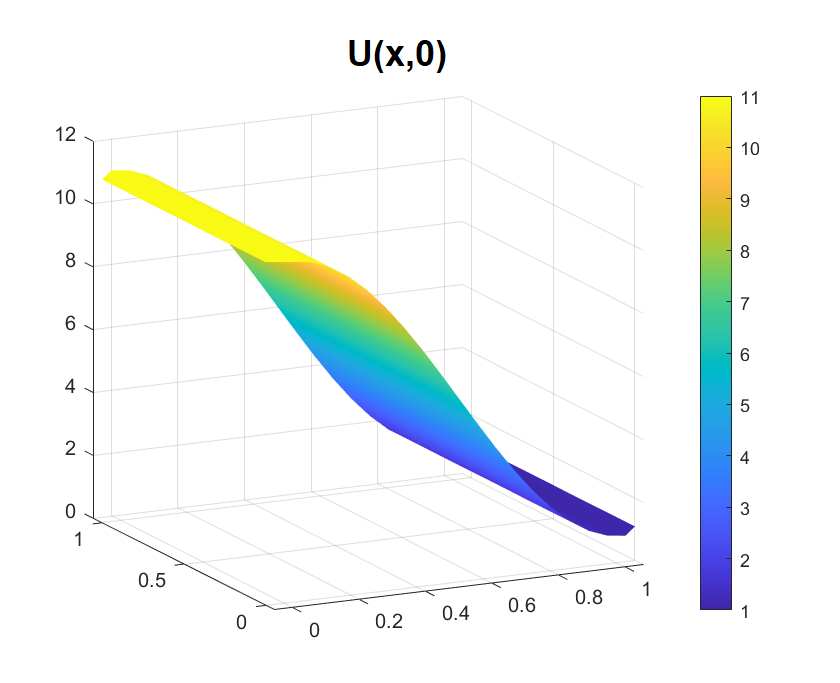}
    \caption{Initial bacterial density from Example 2.}
    \label{fig:5.1}
\end{figure}

\begin{figure}[htp]
    \centering
    \includegraphics[scale=0.3]{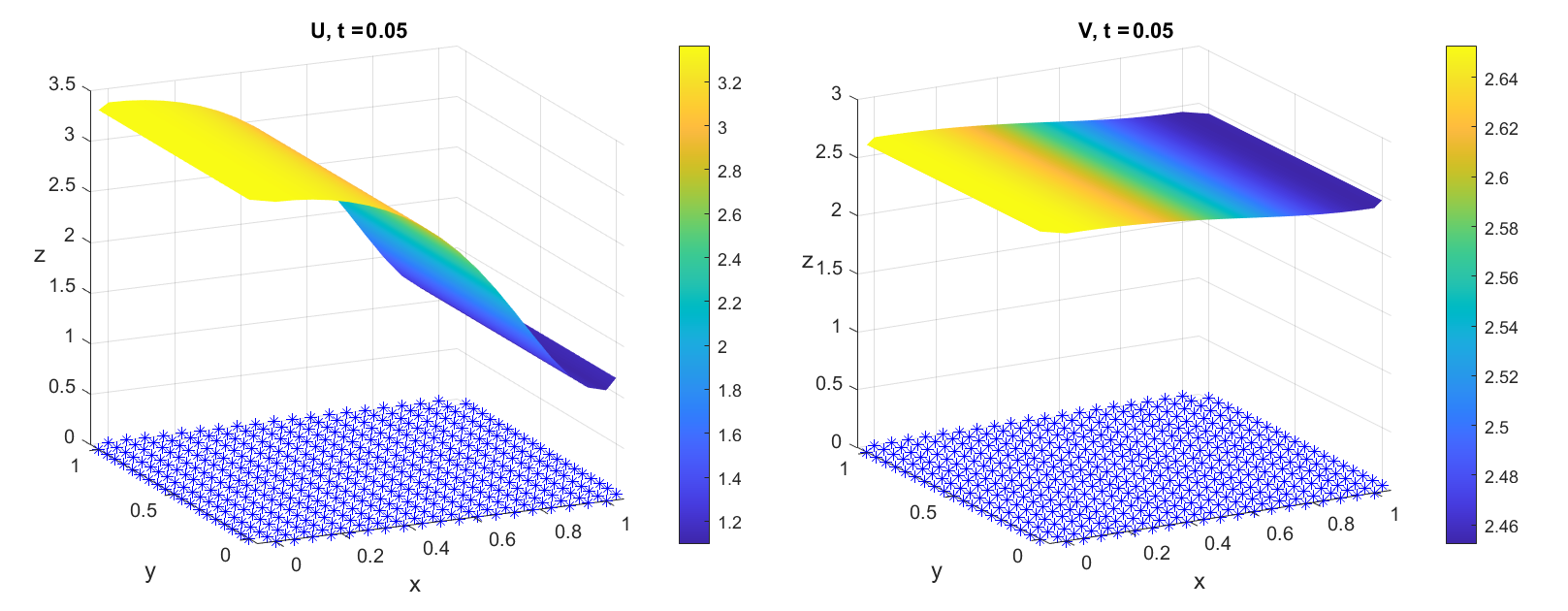}
    \caption{Numerical solution $(U,V)$ from Example 2 at $t = 0.05$ using the regular grid.}
    \label{fig:5.2}
\end{figure}

\begin{figure}[htp]
    \centering
    \includegraphics[scale=0.3]{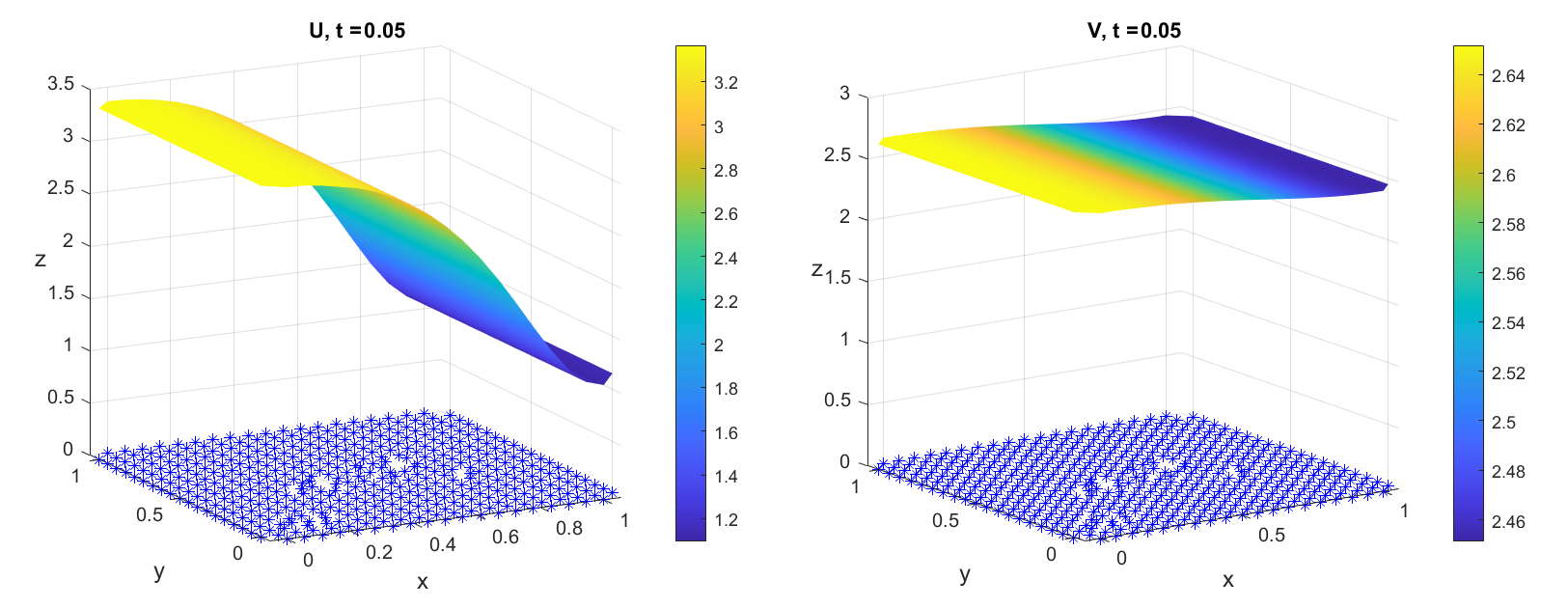}
    \caption{Numerical solution $(U,V)$ from Example 2 at $t = 0.05$ using the irregular cloud of points.}
    \label{fig:5.3}
\end{figure}

We consider again $\Omega = [0,1] \times [0,1]$, this time with initial value $u(\vect{x},0)$ given by:
\begin{equation}\label{5.2}
    u(\vect{x},0) = 6+5\cos(\pi x),
\end{equation}
whose graph is plotted in Figure \eqref{fig:5.1}.
We use a different $\gamma$ function for this example, now being $\gamma(v) := \frac{1}{(1+v)^2}$. Thus, in order to select a $\mu$ large enough for $-2 \gamma'(s) + \gamma''(s)s < \mu$ to hold for all $s \geq 0$, we can take $\mu = 5$, as in this case $-2 \gamma'(s) + \gamma''(s)s = \frac{4}{(1+s)^3} + \frac{6s}{(1+s)^4}$, which is a decreasing function for $s \geq 0$ and takes the value of $4$ at $s=0$. The weights are set as in Example 1 and again we choose $\Delta t = 0.001$.
\\\\
This way, we plot again the results for $t = 0.005$ using the same previous discretizations of $\Omega$. The results for the regular mesh are shown in Figure \ref{fig:5.2}, whereas the irregular cloud of points yields the values of Figure \ref{fig:5.3}.

Lastly, the values of $||U-1||_{l^\infty(\Omega)}$ and $||V-1||_{l^\infty(\Omega)}$ for this example are included in Table \ref{t5.1} at different time instants.
\\\\
As we can check from the values of Table \ref{t5.1}, though the initial value $u(\vect{x},0)$ in this example has a greater $L^\infty$ norm than the one from Example 1, the convergence is faster. Moreover, in this case $u(\vect{x},0) \geq 1$ for all $\vect{x} \in \Omega$, so the initial bacterial distribution is above the carrying capacity of the logistic model (which is 1 after rescaling the variables) at every point. Therefore, the opposite of Example 1 happens, with $u$ converging to $1$ decreasingly.

\begin{table}[h]
\centering
\begin{tabular}{llllll}\hline \\[-1em]
$T(s)$  & 0.05 & 1 & 2.5 & 5 & 10 \\ \\[-1em] \hline
\\[-0.5em]
$||U-1||_{l^\infty(\Omega)}$   & 2.3649 & 0.0051 & $2.6379\cdot 10^{-6}$  & $9.5495 \cdot 10^{-12}$ & $5.8398 \cdot 10^{-14}$  \\ \\[-0.5em]
$||V-1||_{l^\infty(\Omega)}$    & 1.6528 & 0.0049 & $2.6465\cdot 10^{-6}$  & $9.8872 \cdot 10^{-12}$ & $2.3967 \cdot 10^{-12}$  \\ \\[-0.8em] \hline
\end{tabular}
\caption{\label{t5.1} Values of $||U-1||_{l^\infty(\Omega)}$ and $||V-1||_{l^\infty(\Omega)}$ in different time instants from Example 2, calculated with the regular grid.}
\end{table}
\subsection{Example 3}
\begin{figure}[htp]
    \centering
    \includegraphics[scale=0.3]{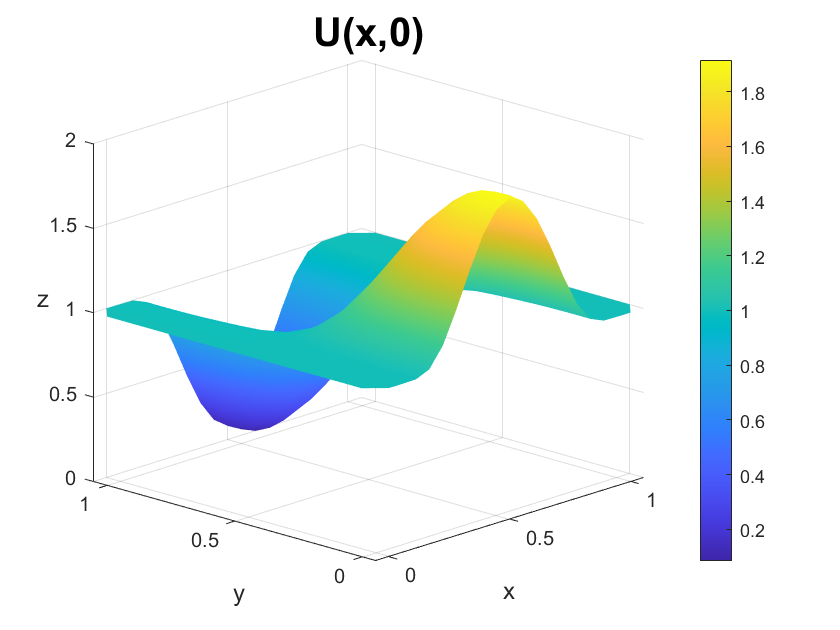}
    \caption{Initial bacterial density from Example 3.}
    \label{fig:5.4}
\end{figure}

\begin{figure}[htp]
    \centering
    \includegraphics[scale=0.4]{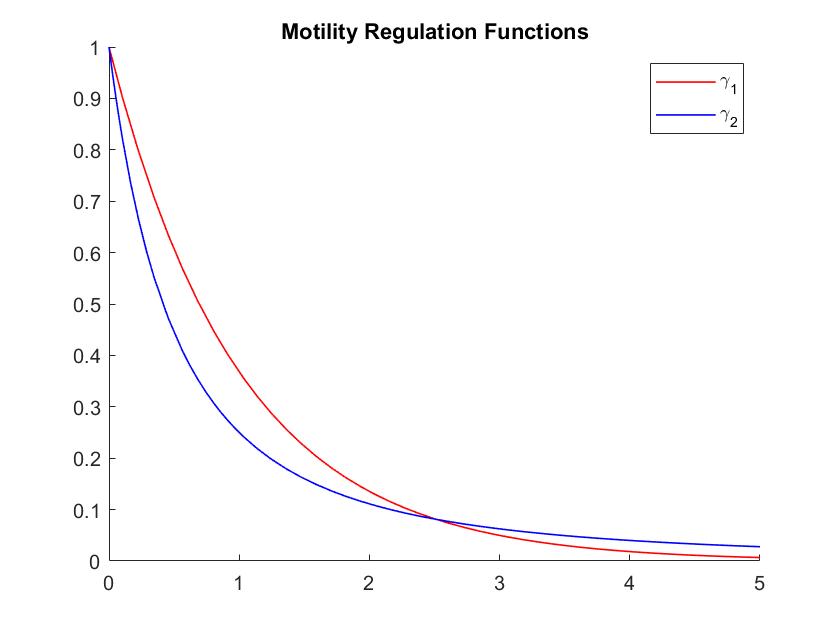}
    \caption{Plots of $\gamma_1(s)$ and $\gamma_2(s)$ for $s \in [0,5]$.}
    \label{fig:5.5}
\end{figure}

For this final example, in order to show a different behaviour than on both previous cases, we choose an initial bacterial distribution $u(\vect{x},0)$ such that $0 <\min_{\vect{x} \in \Omega} u(\vect{x},0) < 1$ and $1 <\max_{\vect{x} \in \Omega} u(\vect{x},0) < \infty$, meaning that for some points in $\Omega$ the initial bacterial density is below the carrying capacity of the logistic model, while it is above for others. As an example, we take:
\begin{equation}\label{5.3}
    u(\vect{x},0) = 1 + 50 \cos(\pi y) e^{-\frac{1}{x(1-x)}} \chi_{(0,1)}(x),
\end{equation}
where
$$
 \chi_{(0,1)}(x) = \begin{cases}
     1, \quad \text{if } x \in (0,1), \\
     0, \quad \text{otherwise},
 \end{cases}
$$
whose graph is shown in Figure \ref{fig:5.4}.
\\\\
This way, the maximum of $ u(\vect{x},0)$ is reached on $(1/2,0)$, with value $1+50 e^{-4} \approx 1.9158$, whereas the minimum lies on $(1/2,1)$, with value $1-50e^{-4} \approx  0.0842$. Therefore, the convergence to $u \equiv 1$ will either be growing or decreasing, depending on the initial value of $u(\vect{x},0)$ on each point. For this example, we compute the numerical solution to problem \eqref{1.2} using both previous choices of $\gamma$, in order to analyze the differences in the behaviour of the solutions. For a clearer notation, we write $\gamma_1(s) = e^{-s}$ and $\gamma_2(s) = \frac{1}{(1+s)^2}$. As can be seen on Figure \ref{fig:5.5}, for small enough values of $s$, $\gamma_1(s)>\gamma_2(s)$, however, after approximately $s=2.5129$, the opposite happens.
\\\\In terms of the biological model, this means that for small enough values of AHL concentration, motility and chemotaxis effects are stronger when the motility regulation is represented by $\gamma_1$. Hence intuitively, with identical initial conditions and value of $\mu$, the solution to problem \eqref{1.2} will more rapidly converge to the homogeneous state $(1,1)$ with $\gamma = \gamma_1$ than with $\gamma = \gamma_2$. To test this numerically we take $\mu = 5$, that fulfills \eqref{1.4} for both choices of $\gamma$, and compute both numerical solutions with the GFD scheme \eqref{4.13}-\eqref{4.14}, using the initial bacterial density given in \eqref{5.3}. The same $\Omega$, $\Delta t$ and $w_i$ as in Examples 1 and 2 are taken.
\\\\
Values of $||U-1||_{l^\infty(\Omega)}$ and $||V-1||_{l^\infty(\Omega)}$ in different time instants are represented in Tables \ref{t5.2} and \ref{t5.3}, for $\gamma = \gamma_1$ and $\gamma = \gamma_2$, respectively.

\begin{table}[h]
\centering
\begin{tabular}{lllllll}\hline \\[-1em]
$T(s)$  & 0.05 & 0.1 & 0.25 & 0.5 & 1 & 2.5 \\ \\[-1em] \hline
\\[-0.5em]
$||U-1||_{l^\infty(\Omega)}$   & 0.4314 & 0.2348 & 0.0577 & 0.0086 & $3.2506 \cdot 10^{-4}$ & $1.2843 \cdot 10^{-7}$ \\ \\[-0.5em]
$||V-1||_{l^\infty(\Omega)}$    & 0.0395 & 0.0315 & 0.0139 & 0.0034 & $2.4541 \cdot 10^{-4}$ & $1.2877 \cdot 10^{-7}$  \\ \\[-0.8em] \hline
\end{tabular}
\caption{\label{t5.2} Values of $||U-1||_{l^\infty(\Omega)}$ and $||V-1||_{l^\infty(\Omega)}$ in different time instants from Example 3, calculated with $\gamma = \gamma_1$ and the regular grid from Figure \ref{fig:4.2}.}
\end{table}

\begin{figure}[htp]
    \centering
    \includegraphics[scale=0.3]{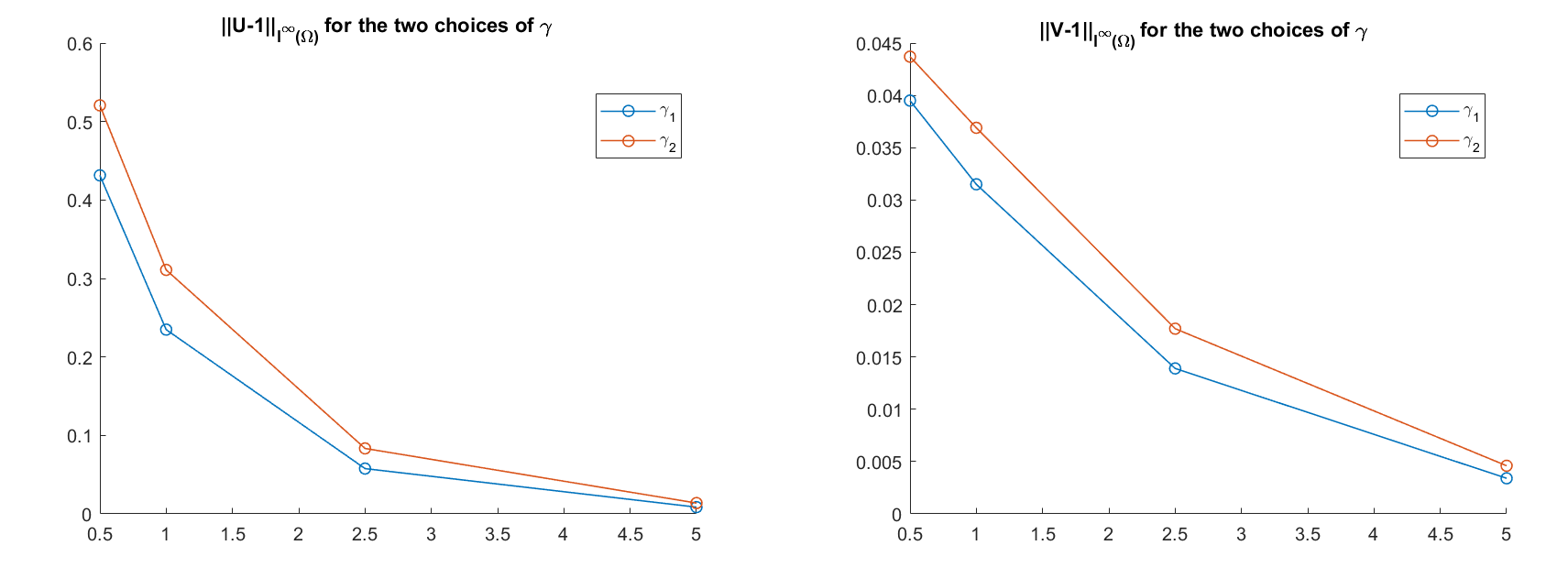}
    \caption{Values of $||U-1||_{l^\infty(\Omega)}$ and $||V-1||_{l^\infty(\Omega)}$ for the two cases considered in Example 3: with $\gamma = \gamma_1$ and $\gamma = \gamma_2$.}
    \label{fig:5.6}
\end{figure}

\begin{table}[h]
\centering
\begin{tabular}{lllllll}\hline \\[-1em]
$T(s)$  & 0.05 & 0.1 & 0.25 & 0.5 & 1 & 2.5 \\ \\[-1em] \hline
\\[-0.5em]
$||U-1||_{l^\infty(\Omega)}$   & 0.5206 & 0.3109 & 0.0834 & 0.0138 & $5.6658 \cdot 10^{-4}$ & $1.6930 \cdot 10^{-7}$\\ \\[-0.5em]
$||V-1||_{l^\infty(\Omega)}$    & 0.0437 & 0.0369 & 0.0177 & 0.0046 & $3.2877  \cdot 10^{-4}$ & $1.6567 \cdot 10^{-7}$ \\ \\[-0.8em] \hline
\end{tabular}
\caption{\label{t5.3} Values of $||U-1||_{l^\infty(\Omega)}$ and $||V-1||_{l^\infty(\Omega)}$ in different time instants from Example 3, calculated with $\gamma = \gamma_2$ and the regular grid from Figure \ref{fig:4.2}.}
\end{table}
By comparing both tables we indeed find that the system with $\gamma_1$ has a faster convergence to the state $(1,1)$ than with $\gamma_2$, especially on the initial moments. This is a consequence of $v$, the AHL concentration, always being below the threshold value $2.5129$, implying that the motility and chemotaxis effects have greater intensity when $\gamma = \gamma_1$. To show these rates of convergence, in Figure \ref{fig:5.6} we plot the previous values of $||U-1||_{l^\infty(\Omega)}$ and $||V-1||_{l^\infty(\Omega)}$ together. To keep a unified scale on the graph, the values are only plotted until $t = 0.5$.
As the graph shows, the red values, corresponding to $\gamma_2$, are always above the blue ones, corresponding to $\gamma_1$. 
\section{Conclusions}
In this work, we've shown the efficiency and simplicity of the Generalized Finite Difference (GFD) Method for solving the non-linear parabolic-elliptic system \eqref{1.2}, modelling the chemotactic response of a strain of \textit{E. coli} bacteria with motility regulation in the presence of AHL. To do so, after a brief introduction to the GFD method, the numerical scheme \eqref{4.13}-\eqref{4.14} was derived, proving its convergence in Theorem \ref{num1}.
\\\\Three different examples were then considered in order to show the applicability of the method. In the first two cases, the GFD method was tested over a regular mesh and an irregular cloud of points, using a different motility regulation function, $\gamma$, on each example, yielding the expected results. Tables \ref{t4.1} and \ref{t5.1} numerically verified the convergence \eqref{1.4.1}. Example 1 was chosen such that the initial bacterial density was below the carrying capacity of the logistic model for all points in $\Omega$, with the solution $u$ hence growing towards the homogeneous state 1. On the other hand, the initial bacterial density of Example 2 was always above the carrying capacity, resulting in a population decrease.
\\\\Finally, in Example 3 the initial density was taken with values both above and below the carrying capacity, so that growth and decrease coexisted in different regions of $\Omega$. Moreover, the two previous choices of $\gamma$ were tested, comparing both convergences.

\section*{Acknowledgments}
This work was partially supported by the Next Generation European Union funds, the Spanish Ministry of Labour and Social Economy and the Ministry of Economy, Finance and Employment of the Community of Madrid (Grant 17-UCM-INV) (F.H-H.). This work is also partially supported by the grant FJC2021-046953-I from MICINN (Spain) (A.M.V.).

\end{document}